\documentclass[12pt]{elsarticle}
\usepackage[left=1in,right=1in, top=1.2in,bottom=1.2in]{geometry}
\usepackage{times}
\usepackage{amssymb,amsthm,latexsym,amsmath,epsfig,pgf}
\usepackage{hyperref}
\usepackage{url}
\usepackage{graphicx}

\newtheorem{theorem}{Theorem}[section]
\newtheorem*{theorem A}{Theorem A}
\newtheorem*{theorem B}{N\"olker's Theorem}
\newtheorem{lemma}{Lemma}[section]

\newtheorem{corollary}{Corollary}[section]

\theoremstyle{remark}
\newtheorem{remark}{Remark}[section]
\theoremstyle{remark}

\begin{document}
 
\begin{frontmatter}
\title{Some results on Bivariegated Graphs and Line Graphs}



\author[label1]{Ranjan N. Naik}

\address[label1]{\small Department of Mathematics,\\
Lincoln University, USA\\

\vspace*{2.5ex} 
 {\normalfont rnaik@lincoln.edu}
 }

\begin{abstract}
This paper is on the structures of line graphs and bivariegated graphs and presents solutions of some graph equations involving line graphs and bivariegated graphs. .

\end{abstract}


\end{frontmatter}

\section{Introduction}
2-variegated (bivariegated) graphs are 1-factorable graphs with some additional properties in [1] Bednarek et al. (1973). Further results on 2-variegated graphs and k-variegated graphs, $ k \geq 2 $, can be found in [2-3], and [6-8]. The earliest works on line graphs are from Whitney (1932) and Krausz (1943) [4].  
We consider here finite graphs without loops or multiple edges and follow the terminology of [4]. Let G (p, q) be a graph with vertex set V and edge set E. A graph B with 2n vertices is said to be a 2-variegated graph if its vertex set can be partitioned into two equal parts such that each vertex from one part is adjacent to exactly one vertex from the other part not containing it [1]. The following characterization of a 2-variegated graph is from [3]. A graph B of order 2n is a 2-variegated graph if and only if there exists a set S of n independent edges in B such that no cycle in B contains an odd number of edges from S. We will call such independent edges in a 2-variegated graph B as special edges.

If  $S$ = $ \{e_{1},e_{2},\ldots,e_{n} \} $ is a set of n special edges of a 2-variegated graph B with 2n vertices, the set S makes a 2-variegation of B. The Petersen graph is a 2-variegated graph. More generally, the same is true for any generalized Petersen graph formed by connecting an outer polygon and an inner star with the same number of points; for instance, this applies to the Mobius-Kantor graph and the Desargues graph.

The line graph L(G) of a graph G is a graph whose vertices are the edges of G and two vertices of L(G) are adjacent if and only if the corresponding two edges of G are adjacent in G.  The characterization of a  line graph G due to Krausz [4] is a graph in which no vertex of G lies on more than two cliques (complete graphs) and every edge of G belongs to exactly one clique.

If G is a 2-variegated graph or not, then the graph L(G) need not be a 2-variegated graph. We solve here the two graph equations:  $L(G)$ = $B$ and  $L(B1)$ = $B2$ where B,  B1 and B2 are 2-variegated graphs. To solve these equations, we have the following lemma 1 and lemma 2.



\begin{lemma}
If  a graph G has an edge $e$ = $uv$ with (degree $d(u) \geq 3$ and degree $d(v) \geq 3$ ) or ( degree $d(u) \geq 3$ and degree $d(v)$ = $1$ ) , then L(G) is not a 2-variegated graph.

Proof.  Since $d(u) \geq 3$ and $d(v) \geq 3$, L(G) has at least 2 cliques of size $\geq 3$, the vertex w in L(G) corresponding to the edge uv in G cannot have a special edge through w because the edges at w in L(G) are on a click of size $ \geq 3.$ This implies that  if uv is an edge in G, and L(G) is a 2-variegated graph, then $d(u) < 3$ or $d(v) < 3.$

Since degree $ d(u) \geq 3 $ and degree $d(v) = 1$, the vertex w in L(G) corresponding to the edge uv in G cannot have a special edge through it. It implies, if uv is an edge in G, and L(G) is a 2-variegated graph for some graph G, then $ d(u) < 2$ or $d(v) \neq 1$. This completes the proof.

\end{lemma}
\begin{lemma}
If G has a cycle of length n, $n \equiv 1, 2, 3 \pmod{4} $, 
then L(G) is not a 2-variegated graph.

Proof: Suppose G has a cycle $x_{1} x_{2},x_{2} x_{3},\ldots,x_{n} x_{1}$ of length n, $n \equiv 1, 2, 3 \pmod{4}$. These edges $x_{i} x_{j}$ are the vertices in L(G) and they are adjacent in L(G) in the same order as in G and form a cycle of length n. If $n \equiv 2 \pmod{4}, $ the alternate $\frac {n}{2}$ edges which are odd in number cannot be the special edges of the vertices on the cycle in L(G) . If $n \equiv 1, 3 \pmod{4} $,  the alternate edges on the cycle cannot be chosen as special edges of the vertices on the cycle in L(G). Let uv be the other edge at $x_{i}$ other than the edges on the cycle in G. Let w be the corresponding vertex in L(G). The edges through w cannot be a special edge of any vertex on the cycle in L(G). It implies, if L(G) is a 2-variegated graph, then G has no cycle of length  $n \equiv 1, 2, 3 \pmod{4}.$ This completes the proof.
\end{lemma}

Two paths are disjoint if they do not have common edges. Let $<$xyz$>$ denote the induced subgraph $P_{3}$ a path of length 2 on the vertices x, y and z with degree d(y) of y is 2.

\begin{theorem}
The graph equation $L(G) = B$ where B is a 2-variegated graph has a solution G if G satisfies the following condition 1. 

$1.$	G (p, q) has even number of edges,  and $\frac{q}{2}$ mutually disjoint $<$xyz$>$ paths such that odd number of these paths do not lie on a cycle.

Conversely, if L(G) is a 2-variegated graph with $p = 2n$ vertices, then G has n  mutually disjoint $<$xyz$>$ paths and odd number of these paths do not lie on a cycle.

Proof.  Proof is by mathematical induction on p the number of vertices. 

Base case: For p = 3, G is $P_{3}$, and L(G) is a 2-variegated graph.

Assume that the result is true for a graph H with $p = k > 3$ vertices satisfying the given condition 1.  Consider a graph $G (r, s)$ with $ r \geq (k + 1) $ vertices and s edges satisfying the given condition. Since there are $\frac{s}{2} $ vertices of degree 2 from the $\frac {s}{2} $ $<$xyz$>$ paths in G, it is clear that the remaining vertices of G are just adjacent to these $\frac{s}{2} $ vertices of degree 2. 

Let v be one of the $\frac{s}{2}$ vertices of degree 2. Let vx and vy be the adjacent edges at v. Let $ S_{x}$ be the set of edges incident at x and $ S_{y} $ be the set of edges incident at y. Delete the vertex v and the edges vx and vy incident at v and also the isolated vertices if any. The resulting graph H satisfies the condition 1 of the theorem and the graph L(H) is a 2-variegated line graph. Let $ V(H)$ = $U_{1} \bigcup V_{1}$ be a 2-variegation with the special edges $ \{e_{1},e_{2}, \ldots,e_{z}\} $. Add a new edge xy to L(H) and join the vertex x to the vertices in one of the sets say $U_{1}$ corresponding to the edges from $S_{x}$ and join y to the vertices in  $V_{1}$ corresponding to the edges from $S_{y}$. The resulting graph L(G) is a line graph of G with a 2-variegation $V(G) = \{U_{1}+x\} \bigcup \{V_{1}+y\}$ with the special edges $\{e_{1},e_{2},\ldots,e_{z},$xy$ \}$. 

Conversely, let L(G) a line graph of G be a 2-variegated graph with $p = 2n$ vertices. Since L(G) has 2n vertices and n special edges, these 2n vertices correspond to the 2n edges of G and each special edge of L(G) corresponds to a pair of adjacent edges say $v_{1} v_{2}$ and $v_{2} v_{3}$ of G and such pairs of edges have to be disjoint because the special edges in L(G) are nonadjacent. By lemma 2, vertex $v_{1}$ cannot be adjacent to vertex $v_{3}$ in G. Let $v_{4} v_{5}$, and $v_{5} v_{6}$ be a pair of edges in G corresponding to some other special edge in L(G). By lemma 1, the edge $v_{4} v_{5}$ or $v_{5} v_{6}$ has to be different from the edges $v_{1} v_{2}$ or $v_{2} v_{3}$. All these pairs of edges in G corresponding to the n special edges in L(G) are disjoint paths of length 2 induced on 3 vertices. So, G has precisely 2n edges and n disjoint paths $<$xyz$>$ of length 2, and degree of y is 2. Also, no odd number of these paths lie on a cycle in G; otherwise L(G) cannot be a 2-variegated graph. This completes the proof of theorem 1.
\end{theorem}

Theorem 1.1 identifies a graph G for which its line graph $L(G)$ is a 2-variegated graph. The following corollary answers when a graph $G$ and $L(G)$ both are 2-variegated graphs.

\begin{corollary}
The equation  $L(B1) = B2$  where B1 and B2 are 2-variegated graphs has a solution B1 if  B1 with 2n vertices satisfies the condition 1 of theorem 1.1 and  has a set of n nonadjacent edges such that each is exactly in one of the n  $<$xyz$>$ paths of length 2.

Proof. It can be verified that the n nonadjacent edges of $B1$ are the n special edges of $B1$ where no odd number of them are on a cycle and thus $B1$ is a 2-variegated graph and $L(B1)$ is also a 2-variegated graph. This completes the proof of the corollary 1.1.
\end{corollary}

\begin{remark}
If B1 has 2n vertices and 2n edges, it has precisely one cycle of length $k \equiv 0 \pmod{4} $.
\end{remark}
Let $L^k(G)$ denote the kth iterated line graph of G.

\begin{corollary}
Only if G is a cycle of length $n \equiv 0\pmod{4}$, then $L^k(G)$ is a 2-variegated graph. 

Proof.  
Case 1: If the max degree of G is 2, then G is a cycle or a path.  
If G is a cycle of length $n \equiv 1, 2, 3\pmod{4},$ then L(G) is not a 2-variegated graph. If G has a cycle of length $n \equiv 0\pmod{4}$,  then L(G) is a 2-variegated cycle of length n and hence $L^k(G)$ is also a 2-variegated graph for $k \geq 1.$

If G is a path of length q, and q is odd, then L(G) is not a 2-variegated graph. If q is even and L(G) is a 2-variegated graph, then L(G) has odd number of edges and hence $L^2(G)$ is not a 2-variegated graph.

Case 2: If G has a vertex u of degree $> 2$, then L(G) has a clique of size $> 2$. If  L(G) is 2-variegated or not, $L^2(G)$ is not a 2-variegated graph. This completes the proof.

\end{corollary}

\begin{corollary}
If G is a 2-variegated graph, then $L(G) = G$ if and only if G is a cycle of length $$n \equiv 0\pmod{4}.$$ 

Proof:  If G is a graph, then $L(G) = G$ if and if only G is a cycle [4]. But a cycle of length $n \equiv 1, 2, 3\pmod{4},$ is not a 2-variegated graph. Hence, $L(G) = G$ if and only if G is a cycle of length n, $$n \equiv 0\pmod{4} .$$ This completes the proof.
\end{corollary}

\end{document}